\documentclass[12pt,a4paper,reqno]{amsart}
\usepackage{amsmath,amssymb,amsfonts,amsthm}
\usepackage{hyperref}
\usepackage[mathscr]{eucal}

\author{H.~M.~Khudaverdian}
\author{Th.~Th. Voronov}
\address{School of Mathematics, University of Manchester, Manchester, M60 1QD, United Kingdom}

\email{theodore.voronov@manchester.ac.uk, khudian@manchester.ac.uk}

\title[Generalized symmetric powers]{On generalized symmetric powers and \\ a generalization of
Kolmogorov--Gelfand--Buchstaber--Rees theory}

\newtheorem*{prop}{Proposition}

\newtheorem*{coro}{Corollary}
\theoremstyle{definition}
\newtheorem*{Rem}{Remark}

\def\co{\colon\thinspace}

\newcommand{\cal}{\EuScript}

\renewcommand{\leq}{\leqslant}
\renewcommand{\geq}{\geqslant}

\DeclareMathOperator{\Ber}{Ber}

\DeclareMathOperator{\Mat}{Mat}

\DeclareMathOperator{\tr}{tr} \DeclareMathOperator{\Tr}{Tr}

\DeclareMathOperator{\str}{str}

\renewcommand{\L}{{\cal L}}

\DeclareMathOperator{\ber}{ber}

\DeclareMathOperator{\Sym}{Sym}

\newcommand{\wed}{\wedge}

\newcommand{\RR}{\mathbb R}

\newcommand{\ZZ}{{\mathbb Z}}

\newcommand{\CC}{\mathbb C}
\newcommand{\NN}{{\mathbb N}}

\newcommand{\f}{\mathbf{f}}

\newcommand{\g}{\mathbf{g}}

\def\a{\alpha}

\newcommand{\F}{{\Phi}}

\newcommand{\rh}{{\boldsymbol{\rho}}}

\newcommand{\ps}{{{\psi}}}

\begin{document}

\maketitle

In 1939, Gelfand and Kolmogorov published a
paper~\cite{gelfand-kolmogorov}, where they showed that for a
(compact Hausdorff) topological space $X$,   homomorphisms of the
algebra of continuous functions $C(X)$ to the field of real numbers
are in a one-to-one correspondence with points of $X$. The algebra
$C(X)$ is considered  purely algebraically, without a topology. This
result is less known than its analog that gave birth to the theory
of normed rings. The Gelfand--Kolmogorov theorem may be viewed as a
description of the image of the canonical embedding of $X$ into the
infinite-dimensional linear space $V=A^*$, where $A=C(X)$, by a
system of quadratic equations $\f(1)=1$, $\f(a)^2-\f(a^2)=0$,
indexed by elements of $A$. This aspect was recently emphasized by
Buchstaber and Rees (see~\cite{buchstaber_rees:2004} and references
therein), who showed that there is a natural embedding into $V$ not
only for $X$, but also for all its symmetric powers $\Sym^n X$. To
this end, algebra homomorphisms should be replaced by the so-called
`$n$-homomorphisms', and quadratic equations describing the image,
by certain algebraic equations of higher degree. Buchstaber and
Rees' theory was motivated by their earlier study of Hopf objects
for multi-valued groups. (In the hindsight, the notion of an
$n$-homomorphism of algebras was present, implicitly, in Frobenius's
notion of  higher  group characters.)

In the present note we suggest a further natural generalization of
Buch\-staber--Rees's theory. Namely, for a set or topological space
$X$ there is a functorial object $\Sym^{p|q} X$, $p,q\geq 0$, and
for a commutative algebra with unit $A$, a corresponding algebra
$S^{p|q}A$ (see definitions below). We call them `{generalized
symmetric powers}'. There is a canonical map from $\Sym^{p|q} X$ to
$V=A^*$;  we introduce certain algebraic equations (see below) that
should describe its image, thus extending the statements of
Gelfand--Kolmogorov and Buchstaber--Rees. On the level of algebras,
this  corresponds  to a description of algebra homomorphisms
$S^{p|q}A\to B$  in terms of the new notion of a
`{$p|q$-homomorphism}'. Our work was motivated by the results in
supergeometry in~\cite{tv:ber}, from which comes our main tool, the
`characteristic function' of a linear map of algebras. The methods
that we propose  yield, in particular, a simple direct proof of the
main theorem of Buchstaber and Rees.

Let $A$ and $B$ be commutative associative algebras with unit.
Consider an arbitrary linear map $\f\co A\to B$. Its
\textit{characteristic function} is defined to be $R(\f,a,z)=e^{\f
\ln(1+az)}$, where $a\in A$ and $z$ is a formal parameter. Example:
if $\f$ is an algebra homomorphism, then $R(\f,a,z)=1+\f(a)z$, i.e.,
a linear polynomial. Algebraic properties of the map $\f$ are
reflected in functional properties of $R(\f,a,z)$ w.r.t. the
variable $z$. Another example: let $R(\f,a,z)$ be a polynomial of
degree $n$. This corresponds to the Buchstaber--Rees theory, as we
show below.

We call  a linear map $\f$, a \textit{$p|q$-homomorphism} if
$R(\f,a,z)$ is a rational function that can be represented by the
ratio of  polynomials of degrees $p$ and $q$. Properties of
$p|q$-homomorphisms follow from general properties of $R(\f,a,z)$.
For an arbitrary map $\f$, $R(\f,a,z)$ has the power expansion at
zero $R(\f,a,z)=1+\ps_1(\f,a)z+\ps_2(\f,a)z^2+\ldots $ where
$\ps_k(\f,a)=P_k(s_1,\ldots,s_k)$. Here $s_k=s_k(\f,a)=\f(a^k)$, and
$P_k$ are the classical Newton polynomials giving expression of
elementary symmetric functions via sums of powers.  There is an
important exponential property $R(\f+\g,a,z)=R(\f,a,z)R(\g,a,z)$.
Suppose $R(\f,a,z)$ extends to a genuine function of $z$ regarded,
say, as a complex variable. Consider its behaviour at infinity. By a
formal transformation one can see that $R(\f,a,z)=z^{\f(1)}e^{\f\ln
a}e^{\f \ln(1+a^{-1}z^{-1})}$. In particular, for $a=1$ we have
$R(\f,1,z)=(1+z)^{\f(1)}$. The assumption that $R(\f,a,z)$ has no
essential singularity  implies that $\f(1)=\chi\in \ZZ$ and the
integer $\chi$ is the order of the pole at infinity. Hence we have
the expansion $R(\f,a,z)=\sum_{k\leq \chi} \ps^*_k(\f,a)z^k$ near
infinity where $\ps_k^*(\f,a):=e^{\f\ln a}\ps_{\chi-k}(\f,a^{-1})$.
We denote the leading term of the expansion $e^{\f\ln
a}=:\ber(\f,a)$ and call it, the \textit{$\f$-Berezinian} of $a\in
A$. Our formal transformations and the definition of
$\mathbf{f}$-Berezinian make sense at least for a rational
characteristic function or for elements in the neighborhood of
identity, if the algebras are over   $\CC$ or $\RR$. One can
immediately see that $\f$-Berezinian is multiplicative:
$\ber(\f,a_1a_2)=\ber(\f,a_1)\ber(\f,a_2)$. Note that $a\mapsto
\ber(\f,a)$ is, in general, a partially defined map  $A\to B$. In
the rational case,  $\ber(\f,a)$ is the ratio of polynomials in the
elements  $\f(a^k)$ with   values in  $B$.

Here are the examples to be kept in mind. If $\rh \co A\to \Mat(n,
B)$ is a matrix representation by $n\times n$ matrices and
$\f(a)=\tr \rh (a)$, then $R(\f,a,z)=\det(1+\rh (a)z)$, and the
$\f$-Berezinian of $a$ is $\det\rh (a)$.  If $\rh \co A\to \Mat(p|q,
B)$ is a matrix representation by matrices $p|q\times p|q$ and
$\f(a)=\str \rh  (a)$ is the supertrace, then $R(\f,a,z)=\Ber(1+\rh
(a)z)$. In this case $\ber(\f,a)=\Ber \rh (a)$ is the ordinary
Berezinian.


Multilinear symmetric functions $\F_{k}(\f,a_1,\ldots,a_k)$ of
elements  $a_i\in A$ such that $\F_{k}(\f,a,\ldots,a)=k!\ps_k(\f,a)$
satisfy relations identical with the Frobenius recursion  for higher
group characters (see~\cite{buchstaber_rees:2004}). The examples
above make this connection apparent. For the case of a matrix
representation, $s_k(\f,a)=\tr\rh(a)^k$, $\ps_k(\f,a)=\tr
\Lambda^k\rh(a)$, and $\F_{k}(\f,a_1,\ldots,a_{k})=k!\,
\tr\left(\rh(a_1)\wed \ldots \wed \rh(a_k)\right)$. The `abstract'
case simply reproduces the relations between these traces.

Let us return to the case when $R(\f,a,z)$ is polynomial in $z$. The
Buchstaber--Rees theory can be recovered as follows.  Here
$\f(1)=\chi=n\geq 0$; it gives the degree of $R(\f,a,z)$, so
$\ps_k(\f,a)=0$ for all $k\geq n+1$ and all $a\in A$. This is
equivalent to the equations $\f(1)=n\in \NN$ and
$\F_{n+1}(\f,a_1,\ldots,a_{n+1})=0$ for all $a_i$, which is
precisely the definition of an $n$-homomorphism according to
Buchstaber and Rees~\cite{buchstaber_rees:2004}. Since here
$\ber(\f,a)=\ps_n(\f,a)$ (in particular, is a polynomial function of
$a$), the latter function is multiplicative in $a$, and therefore
its polarization $\frac{1}{n!}\F_n(\f,a_1,\ldots,a_n)$ yields an
algebra homomorphism $S^nA\to B$. This gives a one-to-one
correspondence between $n$-homomorphisms $A\to B$ and algebra
homomorphisms $S^nA\to B$.

Consider a  topological space $X$. Its  \textit{$p|q$-th symmetric
power} $\Sym^{p|q}X$ is defined as the identification space of
$X^{p+q}$ with respect to the action of $S_p\times S_q$ and the
relations
$$(x_1,\ldots,x_{p-1},y,x_{p+1}\ldots,x_{p+q-1},y)\sim
(x_1,\ldots,x_{p-1},z,x_{p+1}\ldots,x_{p+q-1},z)\,.
$$
The algebraic analog of it is the \textit{$p|q$-th  symmetric power}
$S^{p|q}A$ of a commutative associative algebra with unit $A$. We
define $S^{p|q}A$ as the subalgebra $\mu^{-1}\left(S^{p-1}A\otimes
S^{q-1}A\right)$ in $S^pA\otimes S^qA$ where $\mu\co S^pA\otimes
S^qA \to S^{p-1}A\otimes S^{q-1}A\otimes A$ is the multiplication of
the last arguments. Example: for $A=\CC[x]$, the algebra $S^{p|q}A$
is the algebra of all polynomial invariants of $p|q$ by $p|q$
matrices. There is a relation   between algebra homomorphisms
$S^{p|q}A\to B$ and $p|q$-homomorphisms $A\to B$. To each
homomorphism $S^{p|q}A\to B$ canonically corresponds a
$p|q$-homomorphism $A\to B$, and we have managed to establish the
inverse  in special cases. Example. An element
$[x_1,\ldots,x_{p+q}]\in \Sym^{p|q}X$ defines a $p|q$-homomorphism
on $A=C(X)$: $a\mapsto a(x_1)+\ldots +a(x_{p})-\ldots
 -a(x_{p+q})$. In general, a linear combination of algebra
homomorphisms of the form $\sum n_{\a} \f_{\a}$ where $n_{\a}\in
\ZZ$ is a $p|q$-homomorphism with $\chi=\sum n_{\a}$,
$p=\sum\limits_{n_{\a}>0} n_{\a}$, and $q=-\sum\limits_{n_{\a}<0}
n_{\a}$. This follows from the exponential property of
characteristic function.

By using formulas from~\cite{tv:ber}, the condition that $\f\co A\to
B$ is a $p|q$-homomorphism can be expressed as $\f(1)=p-q$ and
$|\ps_k(\f,a),\ldots,\ps_{k+q}(\f,a)|_{q+1}=0$  for $k\geq p-q+1$,
where $|\ps_k(\f,a),\ldots,\ps_{k+q}(\f,a)|_{q+1}$ is a Hankel
determinant. It is a system of polynomial equations for
`coordinates' of the linear map $\f$. (In particular, it should
describe the image of $\Sym^{p|q}X$.)

The notions of the characteristic function  and $\f$-Berezinian,
introduced in this note, are powerful tools for studying algebraic
properties of linear maps generalizing ring homomorphisms.
Applications of our results  may, in particular, include topological
ramified coverings with branching more general than that considered
by L.~Smith and Dold. The possibility of such an application was
pointed out to us by V.~M.~Buchstaber (compare
with~\cite{buchstaber_rees:2006}). We thank him for fruitful
discussions.

\def\cprime{$'$}

{
\small

\section*{Appendix. Some comments and a short proof of the Buchstaber--Rees main theorem}

In this Appendix we show how our approach allows to give a quick
direct proof of the main theorem of Buchstaber and Rees, namely,
that algebra homomorphisms $S^nA\to B$ are in one-to-one
correspondence with $n$-homomorphisms $A\to B$. For this reason, we
are concerned mainly with the case of a polynomial characteristic
function, which corresponds to the Buchstaber--Rees theory. We also
elucidate some general constructions,  trying not to duplicate the
main text. The questions  of the anonymous referee of the first
version of this note prompted us to write this Appendix, and we
would like to thank him for this.

\section{}
Let $A$ and $B$ be two associative, commutative and unital  algebras
over $\CC$ or $\RR$. Consider a linear map $\f\co A\to B$. We say
that an arbitrary function $\phi\co A\to B$ is \textit{polynomial}
(or $\f$-polynomial) if it is given by a polynomial expression in
$\f(a)$, $\f(a^2)$, etc. It is a useful notion.

The \textit{characteristic function} of a linear map $\f$ introduced
in the main text is defined as
\begin{equation*}
    R(\f,a,z)=\exp\left(\f\left(\ln\left(1+az\right)\right)\right)=1+\psi_1(a)z+
\psi_2(a)z^2+\psi_3(a)z^3+\dots
\end{equation*}
as a formal power series. For brevity denote $R(a,z)=R(\f,a,z)$. It
is a function of both $z$ and $a$. The coefficients $\psi_k(a)$ are
polynomial functions of $a$ of degree $k$, $\psi_k(\lambda
a)=\lambda^k \psi_k(a)$. Indeed, by differentiating the definition
w.r.t. $z$ we can see that $\psi_k(a)$ can be obtained by standard
Newton-like recurrent formulae:
\begin{equation*}
   \psi_1(a) =\f(a), \quad
\psi_{k+1}(a) ={1\over
k+1}\left(\f(a)\psi_k(a)-\f(a^2)\psi_{k-1}(a)+\f(a^3)\psi_{k-2}(a)-\dots\right)\,.
\end{equation*}

Note that the characteristic function according to its definition
obeys the relation
\begin{equation}\label{eq.ident1}
    R(a,z)R(a',z')=R(az+a'z'+aa'zz',1)\,,
\end{equation}
or
\begin{equation}\label{eq.ident2}
    R(a,1)R(b,1)=R(c,1) \quad {\rm if}\quad 1+c=(1+a)(1+b)
\end{equation}
(makes sense as formal power series). We shall use it in what
follows.

\section{}

One can make the following formal transformation of the
characteristic function aimed at obtaining its expansion near
infinity. Note that initially $R(a,z)$ is a formal power series,
which can be seen as the Taylor expansion at zero of some function
of $z$  if such a function exists. Assume that it exists and denote
it by the same $R(a,z)$. Then we have
\begin{multline*}
    R(a,z)=e^{\f\ln(1+az)}=e^{\f\ln\left(az(1+a^{-1}z^{-1})\right)}=\\
    e^{\f\ln(z1)+\f\ln a+\f\ln(1+a^{-1}z^{-1})}=e^{(\ln z)\cdot \f(1)}e^{\f\ln
    a}e^{\f\ln(1+a^{-1}z^{-1})}=\\
    z^{\f(1)}\sum_{k\geq 0}e^{\f\ln a}\psi_k(a^{-1})z^{-k}=\\
    e^{\f\ln a}\,z^{\chi}+e^{\f\ln a}\psi_1(a^{-1})\,z^{\chi-1}+e^{\f\ln
    a}\psi_2(a^{-1})\,z^{\chi-2}+\ldots
\end{multline*}
where we denoted $\chi=\f(1)$. Initially $\f(1)\in B$, but the
assumption that there is a Laurent expansion at infinity forces us
to conclude that $\chi$ must be a number in $\ZZ$.  Here we assume
whatever we need, e.g., that $a^{-1}$ exists, and so on. Instead of
describing how this calculation can be justified, we show below how
one can go around it.

\section{}
Suppose now that the formal power series $R(a,z)$ terminates, i.e.,
it is a polynomial function in $z$ for all $a\in A$. \textit{Under a
mild technical assumption that the degree of $R(a,z)$ is uniformly
bounded by some $N\in \NN$, we deduce that $\f(1)=n\in \NN$ and that
$R(a,z)$ is a polynomial of degree $n$, i.e., the degree is at most
$n$ for all $a$ and is exactly $n$ for some $a$.} The arguments
below may be used to replace the formal calculation above.

Indeed, consider $R(1,z)=\exp[\f(\ln(1+z))]=\exp[\chi\ln(1+z)]$,
where    $\chi=\f(1)\in B$. We shall show first that the element
$\chi\in B$ is an integer. We have $\exp[\chi\ln(1+z)]=(1+z)^\chi$
where $(1+z)^\chi$ is considered is a formal power series:
\begin{equation*}
    (1+z)^\chi=1+\chi z+ {\chi(\chi-1)\over 2}+\dots =
   \sum_{k=0}^\infty {\chi (\chi-1)\dots (\chi-k+1)\over k!}z^k\,.
\end{equation*}
But $R(a,z)$ is a polynomial of degree at most $N$ for all $a$.
Hence $\chi (\chi-1)\dots (\chi-k+1)=0$ for all $k>N$. Under another
technical condition (introduced by Buchstaber and Rees),   that the
algebra $B$ is ``connected'' \footnote{That means by the definition
that the equation $b(b-1)(b-2)\ldots(b-k)=0$ in $B$ implies that
$b=j$ for some $j=0,1,\ldots,k$.}, we conclude that $\chi=n$ for
some integer $n$ between $1$ and $N$.

For an arbitrary $a$, we have, from the above
identity~\eqref{eq.ident1},
\begin{equation*}
    R(a,z)=R(za,1)=R(z-1,1)R\left({1\over z}+a-1,1\right)=z^nR\left({1\over z}+a-1,1\right)\,.
\end{equation*}
More explicitly, the expansion at the RHS has the form
\begin{equation}\label{eq.expans}
    z^n\left[1+\f\left({1\over z}+a-1\right)+\psi_2\left({1\over
z}+a-1\right) +\dots+
  \psi_N\left({1\over z}+a-1\right)\right]
\end{equation}
for some $N\in \NN$. The coefficients $\psi_k(a)$ are polynomial
functions of $a$ of degree $k$, in particular $\psi_k(\lambda
a)=\lambda^k\psi_k(a)$. By comparing the expansions at the LHS and
RHS, we conclude  that the degree of $R(a,z)$ in $z$ can be at most
$n$:
\begin{equation*}
    R(a,z)=1+\psi_1(a)z+ \psi_2(a)z^2+\psi_3(a)z^3+\dots
+\psi_n(a)z^n\,,
\end{equation*}
for any $a$.

Since $\psi_k(a)$ are in one-to-one correspondence with the
functions $\Phi_k(a_1,\ldots,a_k)$ appearing in the ``Frobenius
recursion'' of Buchstaber and Rees (up to the factor of $k!$), which
can be recovered from $\psi_k(a)$ by  polarization \footnote{From
the Frobenius recursion formula, it is easy to prove by induction
that the multilinear maps $\Phi_k$ are symmetric, therefore are
defined by the restrictions to the diagonal, and then, again using
induction, deduce that the functions
$\varphi_k(a)=\Phi_k(a,\ldots,a)$ obey Newton type recurrence
relations and  thus can be identified   with $k!\psi_k(a)$.}, this
reproduces the Buchstaber--Rees definition of $n$-homomorphisms.

\begin{Rem} Here is a formula for the polarization of a homogeneous polynomial
of degree $k$ (the restriction of a symmetric $k$-linear function to
the diagonal):
\begin{equation}\label{eq.polar}
    \Phi_k(a_1,a_2,\dots,a_k)=
              \sum_{r=1}^k\sum_{1\leq i_1<i_2<\dots<i_r\leq
k}(-1)^{k+r}
           \psi_k(a_{i_1}+a_{i_2}+\dots+a_{i_r})\,.
\end{equation}
Here $\Phi_k(a,\ldots,a)=k!\psi_k(a)$. For example,
             $$
              \Phi_2(a_1,a_2)=\psi_2(a_1+a_2)-\psi_2(a_1)-\psi_2(a_2)
             $$
and
\begin{multline*}
    \Phi_3(a_1,a_2,a_3)=\psi_3(a_1+a_2+a_3)-\psi_3(a_1+a_2)-\psi_3(a_1+a_3)-\psi_3(a_2+a_3)\\
    +
  \psi_3(a_1)+\psi_3(a_2)+\psi_3(a_3)\,.
\end{multline*}
This should be a ``textbook formula'', however it is difficult to
find a reference for it.
\end{Rem}

\section{}

We have defined the \textit{$\f$-Berezinian} as a map $\ber_{\f}\co
A\to B$ by the formula $\ber_{\f}(a)=\exp \f(\ln a)=R(a-1,1)$ when
it makes sense. (In the main text, the notation is $\ber(\f,a)$.) It
is clearly multiplicative, for
\begin{equation*}
    \exp \f(\ln ab)=\exp\f\left(\ln a+\ln b\right)=\exp\left(\f\ln
a+\f\ln b\right)=\exp \f(\ln a) \exp \f(\ln b)\,.
\end{equation*}
This holds even if $B$ is non-commutative, but $\f$ is a trace.

Let $\f$ be an $n$-homomorphism. Then  the formally defined
$\ber_{\f}(a)=\exp \f(\ln a)$ is a polynomial function of $a$ with
values in $B$:
\begin{multline*}
    \ber_{\f}(a)=\exp \f(\ln
(1+(a-1)))=R(a-1,1)=\\
1+\f(a-1)+\psi_2(a-1)+\dots+\psi_n(a-1)\,,
\end{multline*}
well-defined for all $a$.
\begin{prop} For an $n$-homomorphism $\f$,
\begin{equation*}
    \ber_{\f}(a)=\psi_n(a)\,.
\end{equation*}
\end{prop}
\begin{proof}
Consider  the equality
\begin{multline*}
    1+\psi_1(a)z+ \psi_2(a)z^2+\psi_3(a)z^3+\dots
+\psi_n(a)z^n=\\
z^n\left[1+\f\left({1\over z}+a-1\right)+\psi_2\left({1\over
z}+a-1\right) +\dots+
  \psi_n\left({1\over z}+a-1\right)\right]
\end{multline*}
(same as~\eqref{eq.expans};     we have legitimately set $N=n$).
Collecting all terms of degree $n$ in $z$, we arrive at
\begin{equation*}
    \psi_n(a)=1+\f\left(a-1\right)+\psi_2\left(a-1\right) +\dots+
  \psi_n\left(a-1\right)=\ber_{\f}(a)\,.
\end{equation*}
\end{proof}
\begin{coro} For an  $n$-homomorphism $\f$, the function
$\psi_n(a)$ is multiplicative in $a$.
\end{coro}

(The original proof of this fact by Buchstaber and Rees was based on
rather hard combinatorial  arguments involving at one instance
hypergeometric polynomials.)

\begin{Rem}
The apparatus of characteristic functions allows to obtain easily
many facts. For example, if $\f$ and $\g$ are $n$- and
$m$-homomorphisms $A\to B$, respectively, then the exponential
property of characteristic functions immediately implies that
$\f+\g$ is an $(n+m)$-homomorphism, since its characteristic
function is the product of polynomials of degrees $\leq n$ and $\leq
m$. If $\g$ is an $m$-homomorphism $A\to B$ and $\f$ is an
$n$-homomorphism $B\to C$, then $R(\f\circ
\g,a,z)=e^{\f\g\ln(1+az)}=e^{\f\ln R(\g,a,z)}=\ber_{\f} R(\g,a,z)$.
Since we know that $R(\g,a,z)$ is a polynomial in $z$ of degree at
most $m$, and the $\f$-Berezinian $\ber_{\f} b$ is a polynomial in
$b\in B$ of degree $n$, we conclude that $R(\f\circ \g,a,z)$ has
degree at most $nm$ in $z$, therefore $\f\circ \g$ is an
$nm$-homomorphism.
\end{Rem}

\section{}

Buchstaber and Rees proved  that there is a one-to-one
correspondence between $n$-homomor\-phisms of $A\to B$ and algebra
homomorphisms $S^nA\to B$.

One can obtain this result in our framework as follows.

Let ${\cal H}^n(A,B)$ be the set of all $n$-homomorphisms from the
algebra $A$ to the algebra $B$.  We shall construct two mutually
inverse maps between the spaces ${\cal H}^n(A,B)$ and ${\cal
H}^1(S^nA,B)$, thus establishing their one-to-one correspondence. It
is convenient to introduce an   $n\times n$ matrix $\L(a)$ with
entries in $A^{\otimes n}=A\otimes \ldots \otimes A$, where
\begin{equation*}
    \L(a)={\rm diag}\,\left[a\otimes 1\otimes \dots \otimes 1,\,
                   1\otimes a\otimes 1\otimes\dots\otimes 1,\dots,
                        1\otimes 1\otimes\dots\otimes 1\otimes a
                        \right]\,.
\end{equation*}
The map $a\mapsto \L(a)$ is a matrix representation.

To every  algebra homomorphism $F\in {\cal H}^1(S^nA,B)$ we assign
an $n$-homomorphism $\f_F\in {\cal H}^n(A,B)$, by setting
\begin{equation*}
    \f_F(a)=F(\Tr \L(a))\,.
\end{equation*}
(Note that $\Tr \L(a)=\Delta (a)$ in the Buchstaber--Rees
notations.)

Conversely, to every  $n$-homomorphism $\f\in {\cal H}^n(A,B)$ we
assign an algebra homomorphism $F_{\f}\in {\cal H}^1(S^nA,B)$
defined by the condition that for all $z$
\begin{equation}\label{eq.det}
   F_{\f}\bigl(\det (1+ \L(a)z)\bigr)=R(\f,a,z)\,.
\end{equation}

We shall prove that these maps are well-defined and  that
\begin{equation*}
    \f_{F_{\f}}=\f,\quad   F_{\f_F}=F\,.
\end{equation*}

Let $F$ be a homomorphism from $S^nA$ to $B$. Then evidently $\f_F$
is a linear map. A calculation shows that the characteristic
function of $\f_F$ is the following polynomial:
\begin{equation}\label{eq.charfF}
    R(\f_F,a,z)=F(\det(1+ \L(a)z))\,.
\end{equation}
It is polynomial of degree $n$. Hence $\f$ is an $n$-homomorphism.
To obtain~\eqref{eq.charfF}, we note that $\f_F\ln(1+az)=F(\Tr
\L(\ln(1+az)))=F(\Tr \ln(1+\L(a)z))$, hence by exponentiating we
arrive at $e^{\f_F\ln(1+az)}=e^{F(\Tr \ln(1+\L(a)z))}=F
e^{\Tr\ln(1+\L(a)z)}=F \det (1+\L(a)z)$.

The formula~\eqref{eq.det} requires a little bit more work. It
defines a map $F_{\f}$ on elements of the form $\Tr\L(a)\in S^nA$,
including   $\det \L(a)$ (which is a polynomial in traces). In
particular,
\begin{equation*}
    F_{\f}(\det\L(a))=\ber_{\f}(a)=\psi_n(a)\,.
\end{equation*}
It is the restriction of the linear map $\frac{1}{n!}\tilde\F_n\co
S^nA\to B$ to the elements of the form $\det
\L(a)=a\otimes\dots\otimes a$, where $\tilde\F_n$ corresponds to the
symmetric multilinear map $\F_n$.  Recall that $\frac{1}{n!}\F_n$ is
the polarization of $\psi_n$, as given by~\eqref{eq.polar}. The
elements $\det \L(a)$ linearly span $S^nA$ and we see that the
linear map $F_{\f}$ is multiplicative on them. Hence it is a
homomorphism.

We proved that the maps $\f\mapsto F_{\f}, F\mapsto \f_F$ are
well-defined. It remains to prove that  their composition is
identity. We have
\begin{equation*}
    F_{\f_F}(\Tr \L(a))=\f_F(a)=F(\Tr \L(a))\,.
\end{equation*}
We see that the maps $F_{\f_F}$ and $F$ coincide on the elements
$\Tr \L(a)$. Hence they coincide on the elements  $\det \L(a)$. This
implies that they coincide on all elements of $S^nA$. To compare
$\f_{F_f}$ and $\f$,  look again at the formula
\begin{equation*}
    R(\f_F,a,z)=F(\det(1+\L(a)z)\,.
\end{equation*}
Expanding it in $z$ and comparing the first terms we see that
$F(\Tr\L(a))=\f(a)$, i.e., $\f_{F_{\f}}=\f$. This concludes the
proof.

}

\end{document}